\documentclass[english, reqno]{ourlematema}

\pdfoutput=1
\usepackage{amsmath}
\usepackage{amssymb}
\usepackage{enumitem}
\usepackage{hyperref}
\usepackage{mathrsfs}
\usepackage{nicematrix}
\usepackage{subcaption}
\usepackage{float}
\usepackage{tikz}
\usepackage{graphics}
\usepackage{bm}
\usepackage{ytableau}
\usepackage{anyfontsize}
\usepackage{mathdots}
\usepackage{tikzscale}
\usepackage{graphicx}

\hypersetup{colorlinks=true,  allcolors=blue, linkcolor=blue,  citecolor=red}

\newtheorem{theorem}{Theorem}
\numberwithin{theorem}{section}
\newtheorem{proposition}[theorem]{Proposition}

\theoremstyle{definition}
\newtheorem{definition}[theorem]{Definition}
\newtheorem{problem}[theorem]{Problem}
\newtheorem{remark}[theorem]{Remark}
\newtheorem{conjecture}[theorem]{Conjecture}
\newtheorem{example}[theorem]{Example}

\DeclareMathOperator{\OGr}{OGr}

\DeclareMathOperator{\Res}{\mathrm{Res}}
\DeclareMathOperator{\sign}{\mathrm{sign}}

\newcommand{\PP}{\mathbb{P}}
\newcommand{\RR}{\mathbb{R}}

\newcommand{\CC}{\mathbb{C} }
\newcommand{\ZZ}{\mathbb{Z}}

\renewcommand{\AA}{\mathbb{A}}

\newcommand{\mcD}{\mathcal{D}}

\newcommand{\mcF}{\mathcal{F}}

\newcommand{\mcP}{\mathcal{P}}

\newcommand{\Gr}{\text{Gr}}

\newcommand{\scrM}{\mathscr{M}}

\title{ \vspace{-5cm} The positive Orthogonal Grassmannian}

\author{\small Yassine El Maazouz }
\address{California Institute of Technology\\ \email{maazouz@caltech.edu}}

\author{\small Yelena Mandelshtam}
\address{Institute for Advanced Study \\ \email{yelenam@ias.edu}}


\date{2024/10/15}

\begin{document}

    \begin{abstract}
		The Pl\"ucker positive region $\OGr_+(k,2k)$ of the orthogonal Grassmannian  emerged as the positive geometry behind the ABJM scattering amplitudes. In this paper we initiate the study of the positive orthogonal Grassmannian $\OGr_+(k,n)$ for general values of $k,n$. We determine the boundary structure of the quadric $\OGr_+(1,n)$ in $\PP^{n-1}_{+}$ and show that it is a positive geometry. We show that $\OGr_+(k,2k+1)$ is isomorphic to  $\OGr_+(k+1, 2k+2)$ and connect its combinatorial structure to matchings on $[2k+2]$. Finally, we show that in the case $n>2k+1$, the \emph{positroid cells} of $\Gr_+(k,n)$ do not induce a CW cell decomposition of $\OGr_+(k,n)$.
	\end{abstract}
	                                                                                                   
	\maketitle
    
	\section{Introduction}
	
	Let $n \geq k$ be positive integers and denote by $\Gr(k,n)$ the \emph{Grassmannian} of $k$-dimensional subspaces of $\CC^n$. The \emph{positive Grassmannian} $\Gr_+(k,n)$ is the semialgebraic set in $\Gr(k,n)$ where all Pl\"ucker coordinates are real and nonnegative. The matroid stratification of the Grassmannian \cite{GGMS} induces a natural decomposition of $\Gr_+(k,n)$ into the so-called \emph{positroid cells}. These cells can be indexed by combinatorial objects like \emph{Grassmann necklaces, decorated permutations, plabic graphs and Le diagrams}, see \cite{postnikov06}.
	
	After Postnikov's landmark paper \cite{postnikov06}, the positive Grassmannian became a rich object of research in algebraic combinatorics \cite{PostnikovICM,GKL22,W21}. Its study accelerated, in recent years, largely due to its unexpected and profound connection to Physics, in particular shallow water waves \cite{KodamaWilliams, AG18} and scattering amplitudes in quantum field theory \cite{TheAmplituhedron,ABCGPT,PositiveGeometries}.
	
	\medskip
	
	The object of study in the present article is the \emph{positive orthogonal Grassmannian} $\OGr_+^{\omega}(k,n)$ which is defined as follows.
    
    \begin{definition}\label{def:OGr}
        Let $\omega: \RR^{n} \times \RR^{n} \to \RR$ be a non-degenerate symmetric bilinear form. We denote by $\OGr^{\omega}(k,n)$ the algebraic variety of \emph{isotropic} $k$-dimensional subspaces $V$ of $\CC^n$ with respect to $\omega$ i.e. $\omega(x,y) = 0$ for any $x,y \in V$. The positive orthogonal Grassmannian $\OGr^\omega_+(k,n)$ is the semi-algebraic subset of $\OGr^\omega(k,n)$ where the Pl\"ucker coordinates are all real and have the same sign.
    \end{definition}
    
	In the special case $n = 2k$ and $\omega(x,y) = \sum_{i=1}^{2k} (-1)^{i-1} x_iy_i$, the semialgebraic set $\OGr^\omega_+(k,2k)$ was first studied in the context of ABJM scattering amplitudes in \cite{huang2014abjm} and later connected to the Ising model in \cite{galashin2020ising}. In this paper we initiate the study of $\OGr^{\omega}_+(k,n)$ for general values of $k,n$ with respect to the quadratic form
	\begin{equation}\label{eq:altForm}
		\omega_0(x,y) = x_1y_1 - x_2y_2 + \cdots + (-1)^{n-1}x_ny_n.
	\end{equation}    
    In particular, we aim to find the combinatorics that govern its boundary structure. We note that the choice of the quadratic form $\omega$ is extremely important. For certain quadratic forms the variety $\OGr^{\omega}(k,n)$ has no real points.
	
	\begin{example}~\\
    \vspace{-5mm}
    \begin{enumerate}
        \item Let $\omega(x,y) = x_1y_1 - x_2 y_2 - x_3 y_3 - x_4 y_4$ be the Lorentzian inner product on $\RR^4$. The variety $\OGr^{\omega}(2,4)$ has no real points. To see why suppose that $V$ is a real point in $\OGr^{\omega}(2,4)$ and let $x,y \in \RR^4$ be a basis of $V$. We have $\omega(x,y) = \omega(x,x) = \omega(y,y) = 0$, so $x_1, y_1 \neq 0$ and we may assume that $x_1 = y_1 = 1$. Then $z = x-y$ satisfies $z_1 = 0$ and $\omega(z,z) = -(z_2^2 + z_3^2 + z_4^2) = 0$. So we deduce that $z=0$ and hence $x = y$ which is a contradiction.
		
	\item Let $\omega(x,y) = x_1y_1 + x_2 y_2 - x_3 y_3 - x_4 y_4$. The variety $\OGr(2,4)^{\omega}$ has real points, for example the rowspan of the matrix:
		\begin{equation}\label{eq:matrixTheta}
			M(\theta) = \begin{bmatrix}
				1 & 0 & \cos(\theta) & -\sin(\theta) \\
				0 & 1 & \sin(\theta) & \cos(\theta) 
			\end{bmatrix} \quad \text{for } \theta \in [0,2\pi].
		\end{equation}
		However, the semialgebraic set $\OGr_+^{\omega}(2,4)$ is zero-dimensional. 
		To see why, if $V$ is a point in $\OGr^{\omega}_+(2,4)$ then the first Pl\"ucker coordinate $p_{12}$ of $V$ does not vanish, otherwise $V$ would contain a non-zero vector of the form $(0,0,a,b)$ and such a vector cannot be orthogonal to itself. Now since $p_{12} \neq 0$, the space $V$ is then the row span of some matrix of the form \eqref{eq:matrixTheta} and one can check that the only such matrix with nonnegative minors is the matrix $M({\pi/2})$.
    \end{enumerate}
	\end{example}

	This article is organized as follows. 
        In Section \ref{sec:2} we collect some facts on the geometry of the orthogonal Grassmannian $\OGr^{\omega}(k,n)$. 
        In particular we determine the ideal of quadrics that cut out $\OGr(k,n)$ in $\PP(\wedge^k \CC^n)$, and determine a Gr\"obner basis for this ideal. 
        In Section \ref{sec3} we investigate $\OGr_+^{\omega_0}(1,n)$ with respect to the alternating form \eqref{eq:altForm}. Namely, we describe its face structure and show that it is a positive geometry. 
        Section \ref{sec4} is devoted to $\OGr^{\omega_0}_+(k,2k+1)$. In this section we show that $\OGr^{\omega_0}_+(k,2k+1)$ is isomorphic to $\OGr^{\omega_0}_+(k+1,2k+2)$ and we relate the face structure of $\OGr^{\omega_0}_+(k,2k+1)$ to matchings on $[2k+2]$. In Section \ref{sec5} we initiate the study of $\OGr^{\omega_0}_+(k,n)$ starting with the case $k=2$. Already in this specific case, we show that the positroid cell decomposition of $\Gr_+(2,n)$ is no longer sufficient to induce a CW cell decomposition of $\OGr_+^{\omega_0}(2,n)$.
        
	\section{Commutative algebra and geometry of \texorpdfstring{$\OGr(k,n)$}{}}
	\label{sec:2}
	
	In this section we collect some facts on the algebraic variety $\OGr^{\omega}(k,n)$ over $\CC$. Since all non-degenerate symmetric bilinear forms over $\CC$ are isomorphic to
	$$(x,y) :=  x_1y_1 + \dots + x_n y_n,$$
	up to a linear change of variables, the varieties $\OGr^{\omega}(k,n)$ for different $\omega$ are isomorphic. So, in this section we may assume that $\omega$ is the standard inner product $(\cdot, \cdot)$, and we suppress $\omega$ and write $\OGr(k,n)$. We recall that the Grassmannian $\Gr(k,n)$ can be embedded in Pl\"ucker space $\PP(\wedge^k \CC^n)$. The $\binom{n}{k}$ coordinates of this projective space are called Pl\"ucker coordinates which we denote by $p_I$ for any subset $I = \{i_1< i_2 < \dots < i_k\}$ of $[n]$. In Physics, the Pl\"ucker coordinates are often denoted by $p_I = \langle i_1 i_2 \dots i_k \rangle$ and for any permutation $\sigma$ of $[k]$ we have $\langle i_{\sigma(1)} i_{\sigma(2)} \dots i_{\sigma(k)}\rangle = \sign(\sigma) \langle i_1i_2\dotsi_k\rangle = \sign(\sigma) p_{I}$.

	\begin{theorem}\label{thm:OGr_Eq}
            The orthogonal Grassmannian $\OGr^\omega(k,n)$ is cut out in $\PP^{\binom{n}{k} - 1}$ by the Pl\"ucker relations in addition to the following $ \frac{1}{2}\binom{n}{k-1} \Big(\binom{n}{k-1}+1\Big)$ equations:
		\begin{equation}\label{eq:OrthogonalEquations}
			\sum_{\ell=1}^{n} \varepsilon(I\ell) \varepsilon(J\ell) \ p_{I\ell} p_{J\ell} = 0, \quad \text{for } I,J \in \binom{[n]}{k-1}.
		\end{equation}
		where $\epsilon(I\ell) = (-1)^{|\{i \in I \colon i > \ell \}|}$ denotes the sign of the permutation that sorts $I\ell$.
	\end{theorem}
	\begin{proof}
		Let $p = (p_I)_{I \in \binom{[n]}{k}}$ be a point in $\OGr(k,n)$ and let $V$ be its corresponding isotropic space in $\CC^n$. The space $V$ is the rowspan of a rank $k$ matrix:
		\[
		\scalebox{0.85}{
			$A = \begin{bmatrix}
				a_{11} & a_{12} & \dots & \dots & a_{1n}\\
				a_{21} & a_{22} & \dots & \dots & a_{2n}\\
				\vdots  & \vdots  & \dots & \dots & \vdots\\
				a_{k1} & a_{k2} & \dots & \dots & a_{kn}
			\end{bmatrix}.$
		}
		\]
		For $I = \{i_1 < \dots <i_k \}$, a Pl\"ucker coordinate $p_I$ can then be written as follows:
		\[
		p_I = \sum_{\sigma \in S_k} \epsilon(\sigma) \prod_{s=1}^{k} a_{\ell j_{\sigma(s)}},
		\]
		where $S_k$ is the symmetric group on $[k]$. Now let $I = (i_1 < \dots <i_{k-1})$ and $J = (j_1 < \dots <j_{k-1})$ be ordered subsets of $[n]$, then:
        \begin{equation}
            \resizebox{.85\hsize}{!}{$\begin{aligned}
			\sum_{\ell=1}^{n} \varepsilon(I\ell) \varepsilon(J\ell) p_{I\ell}p_{J\ell} 
			& = \sum_{\ell = 1}^{n} \sum_{\sigma, \tau \in S_k} \sign(\sigma \tau) \prod_{s = 1}^{k-1} a_{\sigma(s) i_s} \ a_{\tau(s) j_{s}} \ a_{\sigma(k) \ell} \ a_{\tau(k) \ell} \\
			& =  \sum_{\sigma, \tau \in S_k} \sign(\sigma \tau) \prod_{s = 1}^{k-1} a_{\sigma(s) i_s} \ a_{\tau(s) j_{s}}  \underbrace{\left( \sum_{\ell=1}^{n} a_{\sigma(k) \ell} \ a_{\tau(k) \ell} \right)}_{=0}.
		\end{aligned}$}
        \end{equation}
		Since $V$ is isotropic, the last sum in the right-hand-side is $0$ so we deduce that:
		\[
		\sum_{\ell = 1}^{n} \epsilon(I\ell) \epsilon(J\ell) p_{I \ell} p_{J \ell} = 0.
		\]
		Conversely, let $p$ be a point in $\Gr(k,n)$ such that:
		\[
		\sum_{\ell = 1}^{n} \epsilon(I\ell) \epsilon(J\ell) p_{I \ell} p_{J \ell} = 0 \quad \text{for all } I,J \subset \binom{[n]}{k}.
		\]
		Since the quadratic form $(\cdot,\cdot)$ is invariant under the action of the symmetric group $S_n$, we can use the action of this group on $\Gr(k,n)$ and, without loss of generality, assume that $p_{1,2,\dots,k} \neq 0$. So we can write the vector space in $\CC^n$ represented by $p$ as the rowspan of the matrix:
		\[
		\scalebox{0.7}{
			$A = \begin{bmatrix}
				p_{1,2,\dots,k}      &     0    &  \dots  & 0      & (-1)^{k-1} p_{2,3,\dots,k,k+1} & (-1)^{k-1} p_{2,3,\dots,k,k+2} & \dots & (-1)^{k-1} p_{2,3,\dots,k,n}\\
				0      & \ddots   &  \dots  & \vdots & (-1)^{k-2} p_{1,3,\dots,k,k+1} & (-1)^{k-2} p_{1,3,\dots,k,k+2} & \dots & (-1)^{k-2} p_{1,3,\dots,k,n}  \\
				\vdots & \vdots   &  \ddots & \vdots & \vdots & \vdots \\
				0    & \dots    &  0      &    p_{1,2,\dots,k}   & p_{1,2,\dots,k-1,k+1} & p_{1,2,\dots,k-1,k+2} & \dots & p_{1,2,\dots,k-1,n}
			\end{bmatrix}$.
		}
		\]
		From the equations \eqref{eq:OrthogonalEquations}, we can see that the rows of this matrix are orthogonal to themselves and to one another so $p \in \OGr(k,n)$ and this finishes the proof.
	\end{proof}
	
	\begin{example}
		The orthogonal Grassmannian $\OGr(2,5)$ is cut out by a system of $5$ Pl\"ucker relations and 15 relations of the form \eqref{eq:OrthogonalEquations}:
            \begin{equation*}
            \resizebox{\textwidth}{!}{$
            \begin{aligned}
		    &p_{12}^2 + p_{13}^2 + p_{14}^2 + p_{15}^2, \ p_{12}^2 + p_{23}^2 + p_{24}^2 + p_{25}^2, \quad \dots, \quad p_{15}^2 + p_{25}^2 + p_{35}^2 + p_{45}^2\\
			&p_{13} p_{23} + p_{14} p_{24} + p_{15} p_{25}, \ -p_{12} p_{23} + p_{14} p_{34} + p_{15} p_{35}, \ \dots, \ p_{14} p_{15} + p_{24} p_{25} + p_{34} p_{35} 
		\end{aligned}$}
            \end{equation*}
		The ideal generated by these relations is the prime ideal of $\OGr(2,5)$.
	\end{example}
	
	\begin{remark}
        In Physics notation the equations \eqref{eq:OrthogonalEquations} become 
        \[
            \sum_{\ell = 1}^{n} \langle i_1 i_2 \dots i_{k-1} \ell \rangle\langle j_1 j_2 \dots j_{k-1} \ell \rangle = 0.
        \]
        This is because $\epsilon(I\ell) p_{I\ell} = \langle i_1 i_2 \dots i_{k-1} \ell \rangle$. The relations \eqref{eq:OrthogonalEquations} can also be obtained from the so-called \emph{co-circuit matrices} in \cite[Section 3]{EPS24}. We denote by $P$ the $\binom{n}{k-1} \times n$ matrix 
		\[
		P_{I \ell} = \begin{cases} \epsilon(I \ell) \ p_{I \ell}, \quad \text{if } \ell \not \in I \\
			0, \quad \quad \quad \quad \text{otherwise } \ell \in I \end{cases}
		,\quad \text{for} \quad I \in \binom{[n]}{k-1} \text{ and } 1 \leq \ell \leq n.
		\]
		The equations \eqref{eq:OrthogonalEquations} are equivalent to $P P^T = 0$. For a different quadratic form $\omega$ whose matrix in the standard basis of $\CC^n$ is $\Omega$, the equations for $\OGr^{\omega}(k,n)$ are obtained by setting $P \Omega P^T = 0$ instead of \eqref{eq:OrthogonalEquations}. For example for the alternating bilinear form \eqref{eq:altForm}, the equations \eqref{eq:OrthogonalEquations} become:
		\begin{equation}\label{eq:orthRelationsAlt}
			(P \Omega P^T)_{IJ} = \sum_{\ell=1}^{n} (-1)^{\ell - 1} \epsilon(I\ell) \epsilon(J\ell) p_{I\ell}p_{J\ell} = 0 ,\quad \text{for} \quad I, J \in \binom{[n]}{k-1}.
		\end{equation}
	\end{remark}
	
	\begin{proposition}\label{prop:irred} The variety $\OGr(k,n)$ is empty if $n < 2k$. When $n = 2k$ it splits into two irreducible connected components, and it is irreducible when $n > 2k$. Moreover we have:
		\[
		\dim(\OGr(k,n)) = k(n-k) - \binom{k+1}{2} \quad \text{for } n \geq 2k.
		\]
	\end{proposition}
	\begin{proof}
		The proof for the first part of the statement can be found in \cite[Proposition on page 735]{GriffithsHarris}. For the dimension count it is more convenient to work with the quadratic form:
		\[ 
		\omega(x,y) = x_1y_n + x_2 y_{n-1}+ \dots +~x_{n}y_1.
		\]
		Let $\Omega$ be its corresponding matrix, $X = (x_{i,j})$ be a $k\times (n-k)$ matrix of indeterminates, and consider the $k \times n$ matrix $A = [{\rm Id}_k | X] = 0$. The equations $ A \Omega A^{T} = 0$ are equivalent to an expression of each entry in the lower right corner of $X$:
		\[\scalebox{0.9}{
			$X = \begin{bmatrix}
				x_{1,1} & \dots & \dots & \dots  &  x_{1, n-k-1} & {\bm x_{1, n-k}}\\
				x_{2,1} & \dots & \dots & \dots  & {\bm x_{2, n-k-1}} &{\bm x_{2, n-k}}\\
				\vdots & \dots & \dots & \iddots  & \dots & \vdots\\
				x_{k,1} & \dots & {\bm x_{k, n-2k+1}} &  \dots  & \dots & {\bm x_{k, n-k}}\\
			\end{bmatrix}$
		},
		\]
		in terms of all the other entries of the matrix $A$ outside that corner. So we deduce that $\dim(\OGr(k,n)) = \dim(\OGr^{\omega}(k,n)) = k(n-k) - \binom{k+1}{2}$.
	\end{proof}
	
	From now on, we fix positive integers $k,n$ such that $n \geq 2k$. Following \cite{EPS24}, let $Y_{k,n}$ denote \emph{Young's lattice}. This is a poset whose elements are subsets of size $k$ in $[n]$ and the order relation in $Y_{k,n}$ is:
	\[\resizebox{\textwidth}{!}{$5
	\langle i_1 < \dots < i_k \rangle \leq \langle j_1 < \dots < j_k \rangle \quad \text{if} \quad i_1 \leq j_1, \ i_2\leq j_2, \ \dots, \ i_{k-1} \leq j_{k-1} \text{ and } i_k \leq j_k.$}
	\]
	We denote by $\widetilde{Y}_{k,n}$ another copy of Young's lattice. As a set $\widetilde{Y}_{k,n} = \binom{[n]}{n-k}$ and the order relation is given by:
	\[
	[ i'_1 < \dots < i'_{n-k} ] \leq [j'_1 < \dots < j'_{n-k} ] \quad \text{if} \quad i'_1 \geq j'_1, \ \dots, \ i'_{n-k} \geq j'_{n-k}.
	\]
	Finally we denote by $\mcP_{k,n}$ the poset which, as a set, is the disjoint union of $Y_{k,n}$ and $\widetilde{Y}_{k,n}$. All order relations in $Y_{k,n}$ and $\widetilde{Y}_{k,n}$ remain order relations in $\mcP_{k,n}$ and in addition to these relations we have $\binom{2k}{k}$ covering relations:
	\[
	   [j'_1< \dots <j'_{n-k}] \ \ <  \ \ \langle i_1< \dots <i_k\rangle 
	\]
	whenever $\{1,2,3, \dots, 2k\} = \{i_1, \dots, i_k\} \sqcup \{j_1, \dots, j_k\}$ is a partition where the set $\{j_1, \dots ,j_k\}$ is the complement $[n] \setminus \{j'_1, \dots ,j'_{n-k}\}$.
    
    \begin{figure}[ht]
        \centering
        \scalebox{0.6}{
         \begin{tikzpicture}
            
            \draw (2,4) node  {$\langle 56 \rangle$};
            \draw (1,3) node  {$\langle 46 \rangle$};
            \draw (0,2) node  {$\langle 36 \rangle$};
            \draw (-1,1) node {$\langle 26 \rangle$};
            \draw (-2,0) node {$\langle 16 \rangle$};
            
            \draw (1.25,3.25)--(1.75,3.75);
            \draw (0.25,2.25)--(0.75,2.75);
            \draw (-0.75,1.25)--(-0.25,1.75);
            \draw (-1.75,0.25)--(-1.25,0.75);

            \draw (1.25,2.75)--(1.75,2.25);
            \draw (0.25,1.75)--(0.75,1.25);
            \draw (-0.75,0.75)--(-0.25,0.25);
            \draw (-1.75,-0.25)--(-1.25,-0.75);
    
            \draw (2,2) node {$\langle 45 \rangle$};
            \draw (1,1) node {$\langle 35 \rangle$};
            \draw (0,0) node {$\langle 25 \rangle$};
            \draw (-1,-1) node {$\langle 15 \rangle$};

            \draw (1.25,1.25)--(1.75,1.75);
            \draw (0.25,0.25)--(0.75,0.75);
            \draw (-0.75,-0.75)--(-0.25,-0.25);

            \draw (1.25,0.75)--(1.75,0.25);
            \draw (0.25,-0.25)--(0.75,-0.75);
            \draw (-0.75,-1.25)--(-0.25,-1.75);

            \draw (2,0) node {$\langle 34 \rangle$};
            \draw (1,-1) node {$\langle 24 \rangle$};
            \draw (0,-2) node {$\langle 14 \rangle$};

            \draw (1.25,-0.75)--(1.75,-0.25);
            \draw (0.75,-1.25)--(0.25,-1.75);

            \draw (1.25,-1.25)--(1.75,-1.75);
            \draw (0.25,-2.25)--(0.75,-2.75);

            \draw (2,-2) node {$\langle 23 \rangle$};
            \draw (1,-3) node {$\langle 13 \rangle$};

            \draw (1.75,-2.25)--(1.25,-2.75);
    
            \draw (1.25,-3.25)--(1.75,-3.75);
             
            \draw (2,-4) node {$\langle 12 \rangle$};

            \begin{scope}[shift={(7.2,-4)}, rotate=180]

            \draw (2,4) node  {[1234]};
            \draw (1,3) node  {[1235]};
            \draw (0,2) node  {[1245]};
            \draw (-1,1) node {[1345]};
            \draw (-2,0) node {[2345]};

            \draw (1.25,3.25)--(1.75,3.75);
            \draw (0.25,2.25)--(0.75,2.75);
            \draw (-0.75,1.25)--(-0.25,1.75);
            \draw (-1.75,0.25)--(-1.25,0.75);

            \draw (1.25,2.75)--(1.75,2.25);
            \draw (0.25,1.75)--(0.75,1.25);
            \draw (-0.75,0.75)--(-0.25,0.25);
            \draw (-1.75,-0.25)--(-1.25,-0.75);

            \draw (2,2) node {[1236]};
            \draw (1,1) node {[1246]};
            \draw (0,0) node {[1346]};
            \draw (-1,-1) node {[2346]};

            \draw (1.25,1.25)--(1.75,1.75);
            \draw (0.25,0.25)--(0.75,0.75);
            \draw (-0.75,-0.75)--(-0.25,-0.25);

            \draw (1.25,0.75)--(1.75,0.25);
            \draw (0.25,-0.25)--(0.75,-0.75);
            \draw (-0.75,-1.25)--(-0.25,-1.75);

            \draw (2,0) node {[1256]};
            \draw (1,-1) node {[1356]};
            \draw (0,-2) node {[2356]};

            \draw (1.25,-0.75)--(1.75,-0.25);
            \draw (0.75,-1.25)--(0.25,-1.75);

            \draw (1.25,-1.25)--(1.75,-1.75);
            \draw (0.25,-2.25)--(0.75,-2.75);

            \draw (2,-2) node {[1456]};
            \draw (1,-3) node {[2456]};

            \draw (1.75,-2.25)--(1.25,-2.75);
    
            \draw (1.25,-3.25)--(1.75,-3.75);
             
            \draw (2,-4) node {[3456]};
    
            \end{scope}
    
    
            \draw[red] (2.35,0)--(4.65,0);
            \draw[red] (1.35,-1)--(5.65,-1);
    
            \draw[red] (2.35,-2) .. controls (4,-1.5) and (5,-1.5) .. (6.65,-2);
            \draw[red] (0.35,-2) .. controls (2,-2.5) and (3,-2.5) .. (4.65,-2);

            \draw[red] (1.35,-3) -- (5.65,-3);
            \draw[red] (2.35,-4) -- (4.65,-4);
            
        \end{tikzpicture}}
        \caption{The poset $\mcP_{2,6}$ is created from
            $Y_{2,6}$ and $\widetilde{Y}_{2,6}$ by adding the six covering relations in red.}
        \label{fig:P_(2,6)}
    \end{figure}
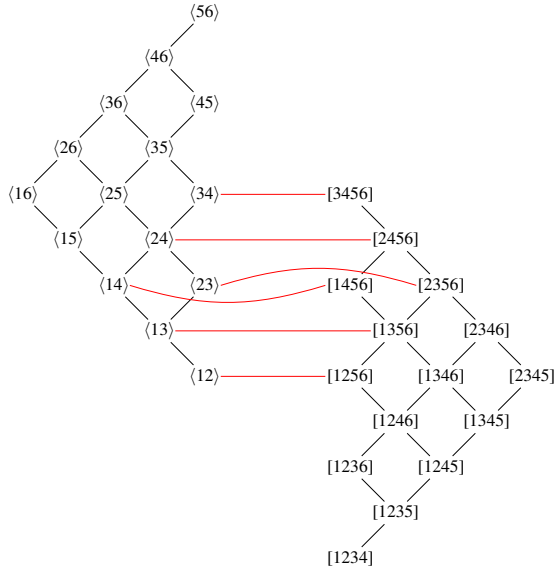
    
	An incomparable pair of elements in $\mcP_{k,n}$ is of type $\big(\langle i_1, \dots, i_k \rangle , \langle j_1, \dots, j_{k}\rangle\big)$ or $\big(\langle i_1, \dots, i_k \rangle , [j'_1, \dots, j'_{n-k}]\big)$. Such a pair yields a non-semistandard Young tableau $\mu$ of shape $(k,k)$ or $\lambda$ of shape $(n-k,k)$:
	
	\setcounter{MaxMatrixCols}{11}
	\begin{equation}\label{eq:skewSymTableau}
		\begin{aligned}
        \mu \ &= 
		\ \begin{bmatrix} 
			j_1 & \cdots &j_{\ell -1}& \bm{j}_{\ell} &\bm{j}_{\ell+1} & \cdots  & \bm{j}_{k}\\
			\bm{i}_1  & \cdots  &\bm{i}_{\ell-1}& \bm{i}_{\ell}  &i_{\ell+1}  & \cdots & i_{k}  \\
		\end{bmatrix},\\
		\lambda \ &= 
		\ \begin{bmatrix} 
			j_1' & \cdots &j'_{\ell -1}& \bm{j}'_{\ell} &\bm{j}'_{\ell+1} & \cdots  &\bm{j}'_{k} &  \cdots & \bm{j}'_{n-k}\\
			\bm{i}_1  & \cdots  &\bm{i}_{\ell-1}& \bm{i}_{\ell}  &i_{\ell+1}  & \cdots & i_{k} & &  \\
		\end{bmatrix}.
		\end{aligned}
	\end{equation}
    The tableau $\mu$ or $\lambda$ being non-semistandard means that there exists an index $\ell$ in $[k]$ such that:
	\begin{equation}\label{eq:snake}
		i_1 < \dots < i_{\ell} < j_{\ell} < \dots < j_{k} \quad \text{or} \quad i_1 < \dots < i_{\ell} < j'_{\ell} < \dots < j'_{n-k}.
	\end{equation} 
	We pick $\ell$ to be the smallest index with this property. The strictly increasing sequences of integers in \eqref{eq:snake} are highlighted in bold in \eqref{eq:skewSymTableau}. Now consider the permutations $\pi$ of the sequence $i_1 < \dots< i_\ell < j_\ell < \dots < j_k$ which make the first $\ell$ entries and the last $k - \ell + 1$ entries separately increasing, and similarly, the permutations $\sigma$ of the sequence $i_1 < \dots< i_\ell < j'_\ell < \dots < j'_{n-k}$ which make the first $\ell$ entries and the last $n- k - \ell + 1$ entries separately increasing. Such permutations permute the bold entries in the tableaux $\mu$ and $\lambda$ in \eqref{eq:skewSymTableau} and yield
    \begin{equation*}
    \resizebox{\textwidth}{!}{$
		\begin{aligned}
        \pi(\mu) \ &= 
		\ \begin{bmatrix} 
			j_1 & \cdots &j_{\ell -1}& \pi(\bm{j}_{\ell}) & \pi(\bm{j}_{\ell+1}) & \cdots  & \pi(\bm{j}_{k})\\
			\pi(\bm{i}_1)  & \cdots  &\pi(\bm{i}_{\ell-1})& \pi(\bm{i}_{\ell})  &i_{\ell+1}  & \cdots & i_{k}  \\
		\end{bmatrix},\\
	\sigma(\lambda) \ &= 
		\ \begin{bmatrix} 
			j_1' & \cdots &j'_{\ell -1}& \pi(\bm{j}'_{\ell}) &\pi(\bm{j}'_{\ell+1}) & \cdots  &\pi(\bm{j}'_{k}) &  \cdots & \pi(\bm{j}'_{n-k})\\
			\pi(\bm{i}_1)  & \cdots  &\pi(\bm{i}_{\ell-1})& \pi(\bm{i}_{\ell})  &i_{\ell+1}  & \cdots & i_{k} & &  \\
		\end{bmatrix}.
		\end{aligned}$}
	\end{equation*}
    Summing over these permutations, the tableaux $\mu$ and $\lambda$ yield quadrics
	\begin{equation} \label{eq:StrLaws} 
        \resizebox{0.93\textwidth}{!}{$
        \begin{aligned}
	    f_{\mu} := \sum\limits_{\pi} \sign(\pi) \ \langle \pi(i_1),\dots, \pi(i_\ell), i_{\ell+1}, \dots i_k\rangle \ \langle j_1,\dots,j_{\ell-1}, \pi(j_{\ell}), \dots, \pi(j_k)\rangle\\
         f_{\lambda} := \sum\limits_{\pi} \sign(\pi) \ \langle \pi(i_1),\dots, \pi(i_\ell), i_{\ell+1}, \dots i_k\rangle \ [j'_1,\dots,j'_{\ell-1}, \pi(j'_{\ell}), \dots, \pi(j'_k)]
	\end{aligned}.
        $}
	\end{equation}
    Here, whenever $J' = \{j'_1<\dots<j'_{n-k}\}$ and $[n]\setminus J' = \{ \bar{j}_1 < \dots < \bar{j}_k \}$ we set
    \[
        [j'_1, \dots, j'_{n-k}] := (-1)^{\sum_{r=1}^{n-k} j'_{r}} \langle \bar{j}_1, \dots, \bar{j}_{k}\rangle.
    \]
    These quadrics are called \emph{straightening laws} in standard monomial theory. We note that, under any reverse lexicographic term order $\prec$ in $\CC[p_{I}]$ given by a linear extension of the poset $Y_{k,n}$, the leading monomial of $f_\mu$ in \eqref{eq:StrLaws} is the binomial $p_{i_1 \dots i_k} \ p_{j_1 \dots j_k}$ and the leading monomial of $f_\lambda$ in \eqref{eq:StrLaws} is the binomial $p_{i_1 \dots i_k} \ p_{\bar{j}_1 \dots \bar{j}_k}$.
	
	\begin{example}\label{ex:YoungTableau} The pair $\langle 12 \rangle$ and $[1356]$ are incomparable in $\mcP_{2,6}$, see Figure \ref{fig:P_(2,6)}. This pair corresponds to the following non-semistandard Young tableau:
		\[
		\lambda = \begin{ytableau} 1 & \textcolor{blue}{\bm{3}}   & \textcolor{blue}{\bm{5}}  & \textcolor{blue}{\bm{6}} \\ 
                  \textcolor{red}{\bm{1}} & \textcolor{red}{\bm{2}} \end{ytableau} \ .
		\]
	    There are $10$ permutations $\pi$ that make the two first entries and the remaining $3$ entries of $\{1,2,3,5,6\}$ separately increasing. Among these $10$ permutations, only $4$ leave $1$ in the bottom row of $\lambda$. The quadric $f_\lambda$ in this case has the following expression:
		\[
		f_\lambda = - \underline{p_{12} p_{24}} - p_{13} p_{24} + p_{15} p_{45} + p_{16} p_{46}.
		\]
	\end{example}
    
	\begin{theorem}\label{thm:Groebner}
		The quadrics in \eqref{eq:StrLaws} form a Gr\"obner basis for the ideal ${\rm I}_{k,n}$ in $\CC[p_{I}]$ generated by the Pl\"ucker relations and the quadratic equations in \eqref{eq:OrthogonalEquations} with respect to any monomial ordering given by a linear extension of the poset $\mcP_{k,n}$.
	\end{theorem}
	\begin{proof}
		The ideal ${\rm I}_{k,n}$ is the image of the ideal $I_{k,n,0}$ of the spinor-helicity variety ${\rm SH}(k,n,0)$ in \cite{EPS24} under the map:
          \[
             \varphi : \CC[p_I, q_{J'}] \mapsto \CC[p_I], \quad   p_I \mapsto p_I, \  q_{J'} \mapsto (-1)^{\sum_{j' \in J'} j'} p_{[n]\setminus J'}.
          \]
        The quadrics \eqref{eq:StrLaws} and \eqref{eq:StrLaws} are the images of a Gr\"obner basis of the ideal $I_{k,n,0}$ \cite[Theorem 2.7]{EPS24} under the map $\varphi$. So we deduce that $f_{\mu}$ and $f_{\lambda}$ in \eqref{eq:StrLaws} and \eqref{eq:StrLaws}  are members of the ideal $I_{k,n}$ for any non-semistandard Young tableaux $\mu$ of type $(k,k)$ and $\lambda$ of type $(n-k,k)$ and the initial monomials $p_{i_1 \dots i_k} \ p_{j_1 \dots j_k}$ and $p_{i_1 \dots i_k} \ p_{\bar{j}_1 \dots \bar{j}_k}$ are in the initial ideal ${\rm in}_{\prec}({\rm I}_{k,n})$. Now let $g \in {\rm I}_{k,n}$ and write $g = \varphi(h)$ for some $h \in I_{k,n,0}$. After scaling $h$ with a suitable nonzero scalar in $\CC$, we can write:
        \[
            h = \underline{\bm{p^a} \bm{q^ b}}  + \sum_{\ell = 1}^{N} c_\ell \bm{p}^{\bm{a}_\ell} \bm{q}^{\bm{b}_\ell},
        \]
        where $\bm{a}, \bm{a_\ell}$ are exponent vectors for the variables $(p_I, \ I \in \binom{[n]}{k})$ and $\bm{b}, \bm{b_\ell}$ are exponent vectors for the variables $(q_{J'}, \ J' \in \binom{[n]}{k})$. The underlined monomial in $h$ is the leading mononmial with respect to $\prec$, so it is divisible by some binomial $p_{I} q_{J'}$ where $I', J'$ is an incomparable pair in $\mcP_{k,n}$. We claim that the leading term of $g = \varphi(h)$ is $\varphi(\bm{p^a} \bm{q^ b})$. To see why let $1 \leq \ell \leq N$, we then have $\bm{p}^{\bm{a}_\ell} \bm{q}^{\bm{b}_\ell} \prec \bm{p^a} \bm{q^ b}$ and after applying the map $\varphi$ this inequality continues to hold.
        
        Hence ${\rm in}_{\prec}(I_{k,n})$ is generated by the binomials $p_{I} p_{J}$ where $I, J' = [n] \setminus J$ is an incomparable pair in $\mcP_{k,n}$. These are exactly the initial monomials of the quadrics \eqref{eq:StrLaws} and \eqref{eq:StrLaws}.   
	\end{proof}

     Our next result gives a formula for the degree of $\OGr(k,n)$ for~$n > 2k$.

    \begin{proposition}\label{prop:degree}
           Let $n > 2k$, $m := \lfloor n/2 \rfloor$, and set $D := k(n-k) - \binom{k+1}{2}$. The degree of $\OGr(k,n)$ in the Pl\"ucker embedding is
           \begin{equation}\label{eq:degreeFormula}
            \resizebox{.915\hsize}{!}{$\begin{aligned}
                  & D! \cdot \left(\prod_{\substack{1 \leq i \leq k \\ k< j \leq m }}  \frac{1}{(2m - i - j)(j-i)}\right)
                \left(\prod_{1 \leq i < j \leq k }  \frac{2}{2m - i - j}\right),  && \text{if } n = 2m, \\
                  & D! \cdot  \left(\prod\limits_{1 \leq i \leq k }  \frac{2}{2m - 2i + 1}\right)  \left(\prod\limits_{\substack{1 \leq i \leq k \\ k< j \leq m }}  \frac{1}{(2m - i - j)(j-i)}\right)
                \left(\prod\limits_{1 \leq i < j \leq k }  \frac{2}{2m - i - j+1}\right),  && \text{if } n = 2m+1.
            \end{aligned}$}
        \end{equation}
        
    \end{proposition}
    \begin{proof}
        
    Fix $n > 2k $ and recall that $\OGr(k,n)$ is irreducible in this case. Let $R$ be homogeneous coordinate ring of $\OGr(k,n)$ endowed with its natural $\ZZ$-grading. By the Borel-Weil-Bott theorem, each degree $\ell$ piece $R_{\ell}$ of $R$ is an irreducible representation $V_{\ell \lambda}$ of ${\rm SO}(n)$ corresponding to the highest weight vector $\ell \lambda$ where $\lambda := e_1 + ... + e_k \in \RR^m$ and $m := \lfloor n/2 \rfloor$.
    Furthermore, we can compute the dimensions of these representations using the Weyl dimension formula as follows. Recall that the positive roots of ${\rm SO}(n)$ are
    \begin{align*}
        \Phi^+ = \begin{cases}
            \{e_i \pm e_j\}_{1 \leq i<j\leq m}, & \text{if } n = 2m\\
            \{e_i \pm e_j\}_{1 \leq i<j\leq m} \cup \{e_i\}_{1\leq i\leq m}, & \text{if } n = 2m+1.
        \end{cases}
    \end{align*}
    The Weyl dimension formula states that for any integer $\ell \geq 0$ we have
    \begin{equation}\label{eq:wdf}
        \dim V_{\ell \lambda} =  \frac{\prod_{\alpha \in \Phi^+} \langle \rho +  \ell \lambda, \alpha \rangle}{\prod_{\alpha \in \Phi^+} \langle \rho, \alpha \rangle} = \prod_{\alpha \in \Phi^+} \left( 1 + \ell \frac{\langle \lambda, \alpha \rangle}{\langle \rho, \alpha \rangle} \right),
    \end{equation}
    where $2\rho$ be the sum of the positive roots $\Phi^+$, and $\langle \cdot, \cdot \rangle$ is the standard inner product on the root space $\RR^m$. This is equal to
    \begin{align*}
        2 \rho = \begin{cases}
             (2m-2)e_1 + (2m-4)e_2 + ... + 2e_{m-1} & n = 2m\\
             (2m-1) e_1 + (2m-3)e_2 + ... + 3e_{m-1} + e_m & n = 2m + 1.
        \end{cases}
    \end{align*}
    Plugging this into the dimension formula above for $n = 2m$ yields
    \begin{equation}\label{eq:dimVl}
        \resizebox{.85\hsize}{!}{$\dim V_{\ell \lambda} = \prod\limits_{\substack{1 \leq i \leq k \\ k < j \leq m}} \left( \left(1 + \frac{1}{2m -i - j} \ell \right)\left(1 + \frac{1}{j -i} \ell \right) \right) \prod\limits_{1 \leq i < j \leq k} \left( 1 + \frac{2}{2m -i -j} \ell \right)$}.
    \end{equation}
    When $n = 2m+1$ we obtain
    \begin{equation*}
        \dim V_{\ell \lambda} = \resizebox{0.85\textwidth}{!}{$
         \prod\limits_{1 \leq i \leq k} \left(1 + \frac{2}{2m-2i+1} \ell \right) \prod\limits_{ \substack{1 \leq i \leq k \\ k < j \leq m} } \left( \left(1 + \frac{1}{2m -i - j + 1}\ell\right)\left(1 + \frac{1}{j -i} \ell \right)\right) \cdot  \prod\limits_{1 \leq i < j \leq k} \left( 1 + \frac{2}{2m -i -j + 1} \ell \right).
        $}
    \end{equation*}    
    This determines the Hilbert polynomial $H(\ell) := \dim(V_{\ell \lambda})$ of the coordinate ring $R$, which is a polynomial in $\ell$ of degree $D$. The degree of $\OGr(k,n)$ is then the leading coefficient of the polynomial $D!  \ H(X)$.
    \end{proof}

    Next we prove that $I_{k,n}$ is a prime ideal for $n > 2k$.

    \begin{lemma}\label{lem:CountPairs}
        The incomparable pairs $(I, J')$ in the poset $P_{k,n}$ with $I \in \binom{[n]}{k}$ and $J' \in \binom{[n]}{n-k}$ such that the pair $(I, [n] \setminus J')$ is comparable in $Y_{k,n}$ are in bijection with semistandard Young Tableaux of shape $(k-1, k-1)$ and fillings in $[n]$. Their number is 
        \[
            \frac{1}{k-1}\binom{n+1}{k}\binom{n}{k-2}.
        \]
    \end{lemma}
    \begin{proof}
        Note that the incomparable pairs $(I, J')$ in question here are exactly the non-semistandard Young tableaux $\lambda$ of shape $(n-k, k)$ with fillings in $[n]$ such that the Young tableau $\lambda^c$ of shape $(k,k)$ obtained by taking the complement of the first row of $\lambda$ in $[n]$ is semistandard. To show the result, it then suffices to exhibit a bijection between the following two sets 
        \[
            \left\{(S_1, S_2): S_1, S_2 \in \binom{[n]}{k-1}, S_1 \leq S_2\right\} 
        \]
         and 
        \[
            \left\{(T_1, T_2): T_1, T_2 \in \binom{[n]}{k}, \ T_1 \leq T_2 \text{ and } T_1^c \not \leq T_2 \right\}.
        \] 
        Here if $A = \{a_1 < \dots < a_\ell \}, B = \{b_1 < \dots < b_\ell\} \in \binom{[n]}{\ell}$, by $A \leq B$ we mean $a_1 \leq b_1, a_2 \leq b_2 ,\dots, a_{\ell} \leq b_{\ell}$. For any $h \in [n]$ and set $S \subseteq [n]$ we denote by $S^h$ the set $S \cap [h]$. Note that $S \leq T$ in Young's lattice if and only if $S^h \geq T^h$ for all $h \in [n]$.
  
         The desired bijection is as follows. Fix $S_1 \leq S_2 \in \binom{[n]}{k-1}$ and let $L = S_1 \cap S_2$ and $R = S_1^c\cap S_2^c$. Let $h$ be minimal such that $|L^h| < |R^h|$. Note that $h$ is guaranteed to exist since $|R^n| = n-2(k-1)+|L^n| > |L^n|$, and further note that we then necessarily have $|L^h| + 1 = |R^h|$. Then we define 
         \[
         T_1 = (S_1 \setminus L^h) \cup R^h, \quad T_2 = (S_2 \setminus L^h) \cup R^h.
         \] 
         It is easy to see that $T_1, T_2 \in \binom{[n]}{k}$ with $T_1 \leq T_2$. Now we show that $T_1^c \not \leq T_2$. Consider the set $(T_1^c)^h = \{t \in T_1^c: t \leq h\} = S_2^h$, while $T_2^h = S_2^h\setminus L^h \cup R^h$. Thus $|(T_1^c)^h| = |S_2^h|$ but $|T_2^h| = |S_2^h| - |L^h|+|R^h| = |S_2^h|+1 > |(T_1^h)^c|$, and therefore $T_1^c \not\leq T_2$.

         The inverse map is as follows. A pair $T_1 \leq T_2 \in \binom{n}{k}$ with $T_1^c\not\leq T_2$ is mapped to the pair
         \[
         S_1 = (T_1 \setminus R^h) \cup L^h, \quad 
         S_2 = (T_2 \setminus R^h) \cup L^h, 
        \]
        where $L = T_1 \cap T_2, \ R = T_1^c\cap T_2^c$, and $h$ is minimal such that $|L^h| > |R^h|$. To see why such an $h$ exists, observe that $T_1^c \not \leq T_2$ implies that there is some $h$ for which $|(T_1^c)^h| < |T_2^h|$, while $T_1 \leq T_2$ implies $|T_2^h| \leq |T_1^h|$. We have 
         \[
         |T_2^h| = |T_2^h \cap T_1^h| + |T_2^h \cap (T_1^c)^h|\quad \text{ and } \quad |(T_1^c)^h| =  |(T_1^c)^h \cap T_2^h|+|(T_1^c)^h \cap (T_2^c)^h|.
         \]
         Together these give $|L^h| = |T_1^h \cap T_2^h| > |(T_1^c)^h \cap (T_2^c)^h| = |R^h|$, as desired. The defined map yields $S_1 \leq S_2 \in \binom{[n]}{k-1}$, so our provided map is a bijection.
    \end{proof}
    \begin{remark}\label{rem:latticepaths}
        We remark that the proof of Lemma \ref{lem:CountPairs} may be more intuitive to the reader if translated into the non-crossing lattice path formulation of semistandard Young tableaux, with a row of a Young tableau defining a path by its left steps \cite[Chapter 7]{sagan}. Under this formulation, taking the complement of a row of a Young tableau corresponds to reflecting the corresponding path. In the proof, the sets $L$ and $R$ correspond respectively to shared left and right steps of the two paths, and the superscript $h$ denotes a cutoff at height $h$. The bijective map reflects the picture of the two lattice paths up to the point at which the number of shared right steps exceeds shared left steps by one (see Figure \ref{fig:latticebijection}).
        
        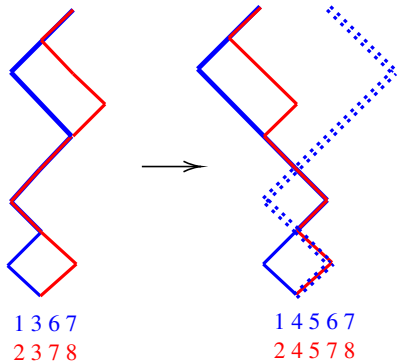
\begin{figure}[ht]
        \begin{center}
        \scalebox{0.8}{
                    \tikzset{every picture/.style={line width=0.75pt}} 
                    
                    \begin{tikzpicture}[x=0.75pt,y=0.75pt,yscale=-1,xscale=1]
                    
                    \draw [color=blue  ,draw opacity=1 ][line width=1.5]    (179.5,191) -- (200,210) ;
                    \draw [color=blue  ,draw opacity=1 ][line width=1.5]    (179.5,191) -- (200.5,170) ;
                    \draw [color=blue  ,draw opacity=1 ][line width=2.25]    (181.5,151) -- (200.5,170) ;
                    \draw [color=blue  ,draw opacity=1 ][line width=2.25]    (181.5,151) -- (199.5,131) ;
                    \draw [color=blue  ,draw opacity=1 ][line width=2.25]    (199.5,131) -- (219.5,110) ;
                    \draw [color=blue  ,draw opacity=1 ][line width=2.25]    (181.5,70) -- (220.5,31) ;
                    \draw [color=red  ,draw opacity=1 ][line width=1.5]    (222.5,190) -- (200,210) ;
                    \draw [color=red  ,draw opacity=1 ][line width=1.5]    (222.5,190) -- (200.5,170) ;
                    \draw [color=red  ,draw opacity=1 ][line width=1.5]    (181.5,151) -- (200.5,170) ;
                    \draw [color=red  ,draw opacity=1 ][line width=1.5]    (199.5,131) -- (181.5,151) ;
                    \draw [color=red  ,draw opacity=1 ][line width=1.5]    (219.5,110) -- (199.5,131) ;
                    \draw [color=red  ,draw opacity=1 ][line width=1.5]    (220.5,71) -- (240.5,90) ;
                    \draw [color=red  ,draw opacity=1 ][line width=1.5]    (200.5,50) -- (220.5,71) ;
                    \draw [color=red  ,draw opacity=1 ][line width=1.5]    (220.5,31) -- (200.5,50) ;
                    \draw [color=blue  ,draw opacity=1 ][line width=1.5]    (339.5,190) -- (360,209) ;
                    \draw [color=blue  ,draw opacity=1 ][line width=1.5]    (339.5,190) -- (360.5,169) ;
                    \draw [color=blue  ,draw opacity=1 ][line width=2.25]    (379.5,150) -- (360.5,169) ;
                    \draw [color=blue  ,draw opacity=1 ][line width=2.25]    (379.5,150) -- (359.5,130) ;
                    \draw [color=blue  ,draw opacity=1 ][line width=2.25]    (359.5,130) -- (298.5,68) ;
                    \draw [color=blue  ,draw opacity=1 ][line width=2.25]    (298.5,68) -- (337.5,29) ;
                    \draw [color=red  ,draw opacity=1 ][line width=1.5]    (382.5,189) -- (360,209) ;
                    \draw [color=red  ,draw opacity=1 ][line width=1.5]    (382.5,189) -- (360.5,169) ;
                    \draw [color=red  ,draw opacity=1 ][line width=1.5]    (379.5,150) -- (360.5,169) ;
                    \draw [color=red  ,draw opacity=1 ][line width=1.5]    (359.5,130) -- (379.5,150) ;
                    \draw [color=red  ,draw opacity=1 ][line width=1.5]    (340,110) -- (359.5,130) ;
                    \draw [color=red  ,draw opacity=1 ][line width=1.5]    (360.5,90) -- (340,110) ;
                    \draw [color=red  ,draw opacity=1 ][line width=1.5]    (318,48.5) -- (360.5,90) ;
                    \draw [color=red  ,draw opacity=1 ][line width=1.5]    (338,29.5) -- (318,48.5) ;
                    \draw [color=blue  ,draw opacity=1 ][line width=2.25]    (181.5,70) -- (219.5,110) ;
                    \draw [color=red  ,draw opacity=1 ][line width=1.5]    (240.5,90) -- (220.5,110) ;
                    \draw    (263.5,129) -- (298.5,129) ;
                    \draw [shift={(300.5,129)}, rotate = 180] [color={rgb, 255:red, 0; green, 0; blue, 0 }  ][line width=0.75]    (10.93,-3.29) .. controls (6.95,-1.4) and (3.31,-0.3) .. (0,0) .. controls (3.31,0.3) and (6.95,1.4) .. (10.93,3.29)   ;
                    \draw [color=blue  ,draw opacity=1 ][line width=1.5]  [dash pattern={on 1.69pt off 2.76pt}]  (358.49,208.89) -- (381.49,187.89)(360.51,211.11) -- (383.51,190.11) ;
                    \draw [color=blue  ,draw opacity=1 ][line width=1.5]  [dash pattern={on 1.69pt off 2.76pt}]  (361.56,167.94) -- (381.56,187.94)(359.44,170.06) -- (379.44,190.06) ;
                    \draw [color=blue  ,draw opacity=1 ][line width=1.5]  [dash pattern={on 1.69pt off 2.76pt}]  (359.47,170.09) -- (339.47,151.09)(361.53,167.91) -- (341.53,148.91) ;
                    \draw [color=blue  ,draw opacity=1 ][line width=1.5]  [dash pattern={on 1.69pt off 2.76pt}]  (339.45,148.93) -- (420.45,68.93)(341.55,151.07) -- (422.55,71.07) ;
                    \draw [color=blue  ,draw opacity=1 ][line width=1.5]  [dash pattern={on 1.69pt off 2.76pt}]  (381.55,28.93) -- (422.55,68.93)(379.45,31.07) -- (420.45,71.07) ;
                    
                    \draw (175,217.4) node [anchor=north west][inner sep=0.75pt]    {$ \begin{array}{l}
                    \textcolor{blue}{1\ 3\ 6\ 7}\\
                    \textcolor{red}{2\ 3\ 7\ 8}
                    \end{array}$};
                    \draw (338,216.4) node [anchor=north west][inner sep=0.75pt]    {$ \begin{array}{l}
                    \textcolor{blue}{1\ 4\ 5\ 6\ 7}\\
                    \textcolor{red}{2\ 4\ 5\ 7\ 8}
                    \end{array}$};
                    
                    \end{tikzpicture}
                    }
        \end{center}
            \caption{A lattice path depiction of the bijection in Lemma \ref{lem:CountPairs}. The red path crosses the reflection of the blue path, pictured as a dotted line.
            \label{fig:latticebijection}}
        \end{figure}
    \end{remark}
    
    \begin{theorem}\label{thm:Primality}
		 When $n > 2k$, the ideal ${\rm I}_{k,n}$ in $\CC[p_{I}]$ generated by the Pl\"ucker relations and the quadratic equations in  \eqref{eq:OrthogonalEquations} is the prime ideal of $\OGr(k,n)$. In particular, the degree of $I_{k,n}$ is given by  \eqref{eq:degreeFormula}.
	\end{theorem}
    \begin{proof}
         By Kostant's theorem, the prime ideal of $\OGr(k,n)$ in $\CC[p_I]$ is generated in degree $2$. To show that ${\rm I}_{k,n}$ is exactly the prime ideal of $\OGr(k,n)$ it is then enough to prove that the dimension of the degree $2$ piece of the $\ZZ$-graded ring $\CC[p_I] / {\rm I}_{k,n}$ is equal to $\dim V_{2 \lambda}$ from the proof of Proposition \ref{prop:degree}. To do so we need to count the dimension of the degree $2$ part of $I_{k,n}$.
                 
         This dimension is equal to the sum of the number of incomparable pairs in the poset $Y_{k,n}$ together with the mixed incomparable pairs in Lemma \ref{lem:CountPairs}. So the dimension of the degree $2$ part of $\CC[p_I] / I_{k,n}$ is 
         \begin{equation}\label{eq:BinomialDim}
             \frac{1}{k}\binom{n+1}{k+1}\binom{n}{k-1} - \frac{1}{k-1}\binom{n+1}{k}\binom{n}{k-2}.
         \end{equation}
         Using \eqref{eq:dimVl}, we compute $\dim(V_{2 \lambda})$ and we find that it is equal to \eqref{eq:BinomialDim}.
    \end{proof}
    \begin{remark}
    \begin{enumerate}
        \item The ideal ${\rm I}_{k,2k}$ is clearly not prime since $\OGr(k,2k)$ has two irreducible connected components and we know that $I_{k,2k}$ cuts out $\OGr(k,2k)$ in $\PP^{\binom{2k}{k}-1}$. Moreover, if $\omega = \omega_0$ is the sign alternating form in \eqref{eq:altForm}, then for any $p \in \Gr(k,2k)$ we have $p \in \OGr^{\omega_0}(k,2k)$ if and only if
        \begin{equation}\label{eq:LinearRelationOGr(k,2k)}
            p_{I} = p_{I^c} \quad \text{for all } I \in \binom{[2k]}{k} \quad \text{ or } p_{I} = - p_{I^c} \quad \text{for all } I \in \binom{[2k]}{k}.
        \end{equation}
        \item Using the \emph{half-spin representation}, $\OGr(k,2k)$ can also be embedded as the spinor variety in $\PP^{2^{n}-1}$. This \emph{spinor embedding} is minimal and the Pl\"ucker embedding we study here can be recovered from the spinor embedding via a quadratic Veronese map. See for example \cite[Section 2]{Manivel} or \cite{Chevalley} for a more detailed account.
    \end{enumerate}
		
    \end{remark}
    
     In this article, we define the \emph{standard} component of $\OGr^{\omega_0}(k,2k)$ to be the connected component where $p_I = p_{I^c}$ for all $I \in \binom{[2k]}{k}$. We will denote the semialgebraic set in the standard component where all Pl\"ucker coordinates are real and have the same sign by $\OGr^{\omega_0}_+(k,2k)$.

	\section{The positive orthogonal Grassmannian \texorpdfstring{$\mathrm{OGr}_{+}(1,n)$}{}}
	\label{sec3}
	
    For the remainder of this article, unless specifically mentioned, we work with the sign alternating quadratic form $\omega_0$ in \eqref{eq:altForm}. We denote by $(p,q)$ the signature of $\omega_0$ where $p = \lceil \frac{n}{2} \rceil$ and $q = \lfloor \frac{n}{2} \rfloor$. In this section we switch gears to study the positive geometry $\OGr_+(1,n) = \OGr^{\omega_0}_+(1,n)$. Here, positive geometry is meant in the sense of \cite[Section 2.1]{PositiveGeometries}. 
    
    We think of the elements of $[n]$ as vertices of a regular $n$-gon ordered clockwise from $1$ to $n$. For each pair of non-empty subsets $A \subset [n]\cap (2\ZZ+1)$ and $B \subset [n] \cap 2\ZZ$, there exists a unique cycle $\sigma(A,B)$ in the symmetric group $S_n$ such that $\sigma(A,B)$ has exactly one excedance and the support of $\sigma(A,B)$ is $A \sqcup B$. The set of such permutations\footnote{These are decorated permutations with all fixed points having a ``$+$" decoration.} $\sigma(A,B)$ is denoted $\mathfrak{S}_{1,n}$. The set $\mathfrak{S}_{1,n}$ is endowed with a partial order given by:
    \[
        \sigma(C,D) \preceq \sigma(A,B) \iff C \subseteq A \text{ and } D \subseteq B.
    \]
    For $\sigma(A,B) \in \mathfrak{S}_{1,n}$, we denote by $\Pi_{\sigma}$ the subset of $\PP^{n-1}_{+}$ where $x_i = 0$ if and only if $i$ is a fixed point of $\sigma(A,B)$ i.e. $i \not \in A \sqcup B$. Here, $\PP^{n-1}_{+}$ is simply $\Gr_+(1,n)$.
    
	\begin{theorem}\label{thm:OGr(1,n)}
		The positive orthogonal Grassmannian $\OGr_{+}(1,n)$ is combinatorially isomorphic to the product of simplices $\Delta_{p-1} \times \Delta_{q-1}$. More precisely, the following hold:

        \begin{enumerate}[leftmargin = 15pt]
            \item $\OGr_+(1,n) = \bigsqcup\limits_{ \sigma \in \mathfrak{S}_{1,n} } \OGr_+(1,n) \cap \Pi_{\sigma}$.
            \item $\overline{\OGr_+(1,n) \cap \Pi_{\sigma}} =  \bigsqcup\limits_{\tau \preceq \sigma} \OGr_+(1,n) \cap \Pi_{\tau}$.
            \item If $A = \{i_1 < \dots < i_r\}$ and $B = \{ j_1 < \dots < j_m\}$ and $\sigma = \sigma(A,B)$ the cell $\OGr_+(1,n) \cap \Pi_{\sigma(A,B)}$ can be parameterized as follows. For each $t_{1},\dots, t_{r-1}$ and $s_{1},\dots, s_{m-1}$ in $\RR_{>0}$ we get a point $x \in \OGr_+(1,n) \cap \Pi_{\sigma(A,B)} $ by setting $x_i = 0$ whenever $i$ is a fixed point of $\sigma(A,B)$ and:
        \end{enumerate}
        \begin{equation}\label{eq:parametrization}
            \resizebox{.915\hsize}{!}{$\begin{aligned}
                x_{i_1} &= \frac{e^{t_1} - e^{-t_1}}{e^{t_1}+e^{-t_1}}, \ &x_{i_2} &= \frac{2}{e^{t_1} + e^{-t_1}} \frac{e^{t_2} - e^{-t_2}}{e^{t_2}+e^{-t_2}}, \ \dots, \ &x_{i_{r-1}} &= \frac{2}{e^{t_{r-1}} + e^{-t_{r-1}}}  \prod\limits_{\ell=1}^{r-1} \frac{e^{t_\ell} - e^{t_\ell}}{e^{t_{\ell}} + e^{-t_{\ell}}},\\
                x_{j_1} &= \frac{e^{s_1} - e^{-s_1}}{e^{s_1}+e^{-s_1}}, \ &x_{j_2} &= \frac{2}{e^{s_1} + e^{-s_1}} \frac{e^{s_2} - e^{-s_2}}{e^{s_2}+e^{-s_2}}, \ \dots, \ &x_{j_{m-1}} &= \frac{2}{e^{s_{m-1}} + e^{-s_{m-1}}}  \prod\limits_{\ell=1}^{m-1} \frac{e^{s_\ell} - e^{s_\ell}}{e^{s_{\ell}} + e^{-s_{\ell}}}.
            \end{aligned}$}
        \end{equation}
        
	\end{theorem}

	\begin{proof}
    The semi-algebraic set $\OGr_{+}(1,n)$ in $\PP^{n-1}$ is cut out by:
      \[
        \sum_{i \in [n]\cap (2\ZZ + 1)} x_i^2  = \sum_{j \in [n]\cap 2\ZZ} x_j^2   \quad \text{and} \quad  x \in \PP^{n-1}_{+}.
      \]
    The cells $\OGr_+(1,n) \cap \Pi_{\sigma(A,B)}$ correspond to the boundaries of $\OGr_+(1,n)$ where the point $x \in \PP_{+}^{n-1}$ satisfies:
    \begin{equation}\label{eq:ABequation}
        \sum_{i \in A} x_i^2 = \sum_{j \in B} x_j^2 \quad \text{ and } \quad x_i \neq 0 \iff i \in A \sqcup B.
    \end{equation}
    The closure of such a cell is obtained by driving some of the coordinates indexed by $A$ or $B$ to $0$. So the first and second statement follow. The face structure described in here is the same face structure as the product $\Delta_{p-1} \times \Delta_{q-1}$ of simplices of dimensions $p-1$ and $q-1$.
    
    For the third statement, to parametrize the cell $\OGr_+(1,n) \cap \Pi_{\sigma(A,B)}$ it is enough to parametrize the set of points $(x_i)_{i \in A} \in \RR_{>0}^{|A|}$ and $(x_j)_{j \in B} \in \RR_{>0}^{|B|}$ on the unit spheres $S^{|A|-1}$ and $S^{|B|-1}$. This is exactly what is done in \eqref{eq:parametrization}.
	\end{proof}

    \begin{example} The orthogonal Grassmannian $\OGr_+(1,5)$ has the same combinatorial structure as $\Delta_{1} \times \Delta_{2}$. The poset of the boundaries of $\OGr_+(1, 5)$ is depicted in Figure \ref{fig:OGr15}.
    \begin{figure}[ht]
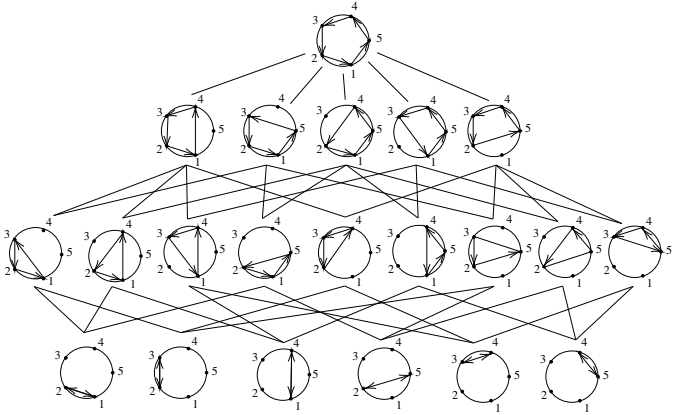

        \centering
        \scalebox{0.5}{

\tikzset{every picture/.style={line width=0.75pt}} 


        }
        \caption{The Hasse diagram of the poset structure on $\mathfrak{S}_{1,5}$.}
        \label{fig:OGr15}
    \end{figure}
    \end{example}

The next theorem shows that $\OGr_+(1,n)$ is a positive geometry. To give its canonical form, it is convenient to permute\footnote{Here, since $k=1$, permuting the coordinates does not change the signs of the ``minors".} the coordinates of $\PP^{n-1}$ and write $\OGr_+(1,n)$ as follows:
\[
    \OGr_+(1,n) = \big\{(x_1: \dots: x_n) \in \PP_{+}^{n-1} \colon x_1^2 + \dots + x_p^2 - x_{p+1}^2 - \dots - x_n^{2} = 0 \big\}.
\]

 \begin{theorem}
         The semi-algebraic set $\OGr_+(1,n) \subset \PP_{+}^{n-1}$ is a positive geometry in the sense of \cite[Section 2.1]{PositiveGeometries}. Its canonical form:
		\[
		\Omega = (1+u_{2,1}^2 + u_{3,1}^2 + \cdots + u_{p,1}^2) \  \frac{ du_{2,1} \wedge du_{3,1} \wedge \dots \wedge du_{n-1,1} }{ u_{2,1} \ u_{3,1} \ \ \cdots \ u_{n-1,1} \ u_{n,1}^2}.
		\]
		Here, $u_{i,j} := x_i/x_1$ in the projective coordinates $(x_1 : \dots : x_n)$ of $\PP^{n-1}$.
    \end{theorem}
    \begin{proof}
    For ease of notation we write $u_i$ for $u_{i,1}$. In the way the form $\Omega$ is defined, it is clear that for $p+1\leq i \leq n-1$:
    \[
        \Res_{u_i = 0} (\Omega) = (-1)^{i+1} (1+u_2^2 + \dots +u_p^2) \frac{du_2 \wedge \dots\wedge \widehat{du_i} \wedge \dots \wedge du_{n-1}}{u_2 \dots \widehat{u_i} \dots u_{n-1} u_n^2}.
    \]
    For $2 \leq i \leq p$ we can rewrite $\Omega$ as:
    $
    \Omega = -(u_{p+1}^2 + \dots + u_n^2) \  \frac{ du_2 \wedge du_3 \wedge \dots \wedge du_{n-1} }{ u_{2} \ u_{3} \ \ \cdots \ u_{n-1} \ u_{n}^2},
    $
    and then:
    \[
        \Res_{u_i = 0} (\Omega) = (-1)^{i+1} (-u_{p+1}^2 - \dots - u_n^2) \frac{du_2 \wedge \dots \wedge \widehat{du_i} \wedge \dots \wedge du_{n-1}}{u_1 \dots \widehat{u_i} \dots u_{n-1} u_n^2}.
    \]
    For the boundary corresponding to $u_n = 0$ note that we have:
    \[
        u_2 du_2 + \dots + u_p du_p - u_{p+1}du_{p+1} - \dots - u_n du_n = 0.
    \]
    Taking the wedge on the right with $du_3 \wedge \dots \wedge du_{n-1}$ we get:
    \[
        u_2 du_2 \wedge \dots \wedge du_{n-1} - u_n du_n \wedge du_3 \wedge \dots \wedge du_{n-1} = 0.
    \]
    So we can rewrite $\Omega$ as follows:
    \[\resizebox{.94\hsize}{!}{$
        \Omega = (1+u_{2}^2 + \cdots + u_{p}^2) \  \frac{ du_{2} \wedge du_{3} \wedge \dots \wedge du_{n-1} }{ u_{2} \ u_{3} \ \ \cdots \ u_{n-1} \ u_{n}^2}
        =  (1+u_{2}^2 + \cdots + u_{p}^2) \  \frac{du_n \wedge du_3 \wedge \dots \wedge du_{n-1} }{u_2^2 u_3 \dots u_{n-1} u_n}$.}
    \]
    We then get the residue at $u_n = 0$:
    \[
        \Res_{u_n = 0} (\Omega) = (1+u_{2}^2 + \cdots + u_{p}^2) \  \frac{ du_3 \wedge \dots \wedge du_{n-1} }{u_2^2 u_3 \dots u_{n-1}}.
    \]
    The residue at the boundary where $x_1 = 0$ can be computed similarly by switching to a different affine chart (eg. set $x_2 = 1$). This covers all the boundaries of $\OGr_+(1,n)$ and since the residue at each boundary gives the same form for a lower dimensional positive orthogonal Grassmannian $\OGr_{+}(1, n)$, we can carry on taking residues for lower dimensional boundaries in the same way.
    \end{proof}

	\begin{example}[$\OGr_{+}(1, 4)$]
    The points $(x_1:x_2:x_3:x_4)$ in the positive orthogonal Grassmannian $\OGr_{+}(1,4)$ in $\PP^3$ are those that satisfy:
		\[
		x_1^2 - x_2^2 + x_3^2 - x_4^2 = 0 \text{ and } x_1,x_2,x_3,x_4 \geq 0.
		\]
		We can see that $\OGr_{+}(1,4)$ is a curvy quadrilateral inside the $3$-simplex in $\PP^3_{+}$ as depicted in Figure \ref{fig:OGr+1,4}. It is a positive geometry and its canonical form is:
		\begin{align*}
			\Omega &= (1 + u_{3,1}^2) \ \frac{du_{2,1} \wedge du_{3,1}}{u_{3,1} \ u_{2,1} \ u_{4,1}^2 }.
		\end{align*}
        where $u_{i,j} = x_i/x_j$.
		
		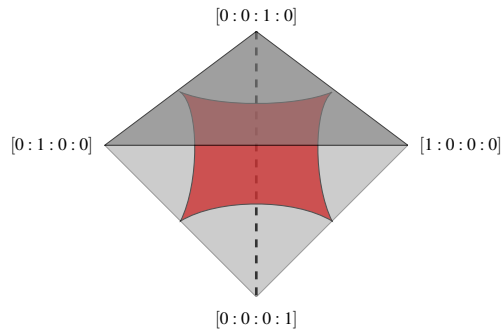
\begin{figure}[ht]
			\centering
				\begin{tikzpicture}
					\coordinate (a) at (-2,0);
					\coordinate (b) at (2,0);
					\coordinate (c) at (0,1.5);
					\coordinate (d) at (0,-2);
					
					\draw[draw=black, line width=1pt, dashed] (c) -- (d);
					
					\draw[black] (2.7,0) node  [xscale = 0.6, yscale = 0.6] {$[1:0:0:0]$};
					\draw[black] (-2.7,0) node  [xscale = 0.6, yscale = 0.6] {$[0:1:0:0]$};
					\draw[black] (0,1.7) node  [xscale = 0.6, yscale = 0.6] {$[0:0:1:0]$};
					\draw[black] (0,-2.3) node  [xscale = 0.6, yscale = 0.6] {$[0:0:0:1]$};

					\filldraw[draw=black, fill=red, opacity = 0.8] 
					(1,0.7) .. controls (0.5,0.5) and (-0.5,0.5) .. (-1,0.7)  
					.. controls (-0.75,0.5) and (-0.75,-0.5) .. (-1,-1)
					.. controls (-0.5,-0.7) and (0.5,-0.7) .. (1,-1)
					.. controls (0.75,-0.5) and (0.75,0.5) .. (1,0.7);
					\filldraw[draw=black, fill=gray, opacity = 0.75, line width=0pt] (a) -- (b) -- (c) -- cycle;
					\draw[draw=black, fill=gray, opacity = 0.4, line width=0pt] (a) -- (b) -- (d) -- cycle;        
				\end{tikzpicture}
			\caption{The positive Grassmannian $\OGr_{+}{(1,4)}$ is the red region in the tetrahedron $\PP^3_{+}$. The boundaries of $\OGr_{+}(1,4)$ lie on the facets of $\PP^3_{+}$.}
			\label{fig:OGr+1,4}
		\end{figure}
		
	\end{example}
 
	\section{The positive Orthogonal Grassmannian \texorpdfstring{$\OGr_+(k,2k+1)$}{}}
	\label{sec4}
    
	We recall that we are working with the sign alternating form \eqref{eq:altForm}. 
    The positroid cells of $\Gr_+(k,2k)$ induce a CW cell decomposition on the positive orthogonal Grassmannian $\OGr_+(k,2k)$, and the cells of this decomposition are indexed by fixed-point-free involutions of $[2k]$. 
    The face structure of $\OGr_+(k,2k)$ as well as the parametrization of its cells are studied in detail in \cite[Section 5]{galashin2020ising}.
    
    One of the reasons positroid cells induce a cell decomposition of the positive orthogonal Grassmannian $\OGr_+(k,2k)$ is that the latter can be obtained by slicing $\Gr_+(k,2k)$ by a linear space, see \eqref{eq:LinearRelationOGr(k,2k)}. In general, one can obtain $\OGr_+(k,n)$ by slicing the positive flag variety with a linear space as we shall now explain.
    
    For a subspace $V$ in $\CC^n$ of dimension $k$, we denote by $V^\perp$ its orthogonal complement with respect to the form \eqref{eq:altForm}.
    
	\begin{lemma}\label{lem:hodgeStarMap}
		The Hodge star map $\Gr(k,n) \to \Gr(n-k,n), V \mapsto V^\perp$ is given in Pl\"ucker coordinates by:
		\[
		q_J = p_{J^c}, \quad \text{for any } J \in \binom{[n]}{n-k},
		\]
		where the $p_I$ and $q_J$'s are the Pl\"ucker coordinates in $\Gr(k,n)$ and $\Gr(n-k,n)$ respectively. In particular it restricts to an isomorphism of positive geometries between $\Gr_+(k,n)$ and $\Gr_+(n-k,n)$.
	\end{lemma}
	\begin{proof}
		If $(p_{I})_{I}$ are the Pl\"ucker coordinates of $V$, the Pl\"ucker coordinates $q_{J}$ of $V^\perp$ are obtained as follows:
		\[
		q_J = s_J \ t_J \ p_{J^c} \quad \text{for any } J \in \binom{[n]}{n-k},
		\]
		where $t = \prod_{j \in J} \omega_0(e_j,e_j)$ and $s_{J}$ is the sign of the permutation that sorts the string $J J^c$. We compute the signs $s_J$ and $t_J$ and get the following:
		\[
		s_J = (-1)^{\sum_{j \not \in J} j  - k(k+1)/2} \quad \text{ and } t_J = (-1)^{|J \cap 2\ZZ|}.
		\]
		From this it is not so hard to see that the sign $s_J t_J$ does not depend on $J$. Since we work with projective coordinates, we deduce that $q_J = p_{J^c}$, hence the Hodge star map is indeed an isomorphism between $\Gr_+(k,n)$ to $\Gr_+(n-k,n)$.
	\end{proof}

    Let $\mcF(k,n)$ be the $2$-step flag variety of partial flags $V \subset W \subset \CC^n$ where $\dim(V) = k$ and $\dim(W) = n-k$. The positive part $\mcF_+(k,n)$ of $\mcF(k,n)$ is the semi-algebraic set of points $(V,W) \in \Gr_+(k,n) \times \Gr_+(n-k, n)$ such that $(V,W) \in \mcF(k,n)$. We denote by $\mcD$ the \emph{diagonal} subset of $\PP^{\binom{n}{k}} \times \PP^{\binom{n}{n-k}}$ i.e.:
	\[
	\mcD := \left\{(p, q) \colon p_I = q_{I^c} \quad \text{for any } I \in \binom{[n]}{k} \right\}.
	\]

    \begin{proposition}
    The positive orthogonal Grassmannian $\OGr_+(k,n)$ is the intersection of the positive flag variety $\mcF_+(k,n)$ with $\mcD$ i.e.:
    	\begin{equation}\label{eq:DiagSlice}
    		\OGr_+(k,n) = \mcF_+(k,n) \cap \mcD.
    	\end{equation}
    \end{proposition}
    \begin{proof}
        This follows immediately from Lemma \ref{lem:hodgeStarMap}.
    \end{proof}
 
	This motivates the choice of the sign alternating form \eqref{eq:altForm} in \cite{galashin2020ising,huang2014abjm}. However, unlike $\mcF_+(k,2k) \cong \Gr_+(k,2k)$, the positive region $\mcF_+(k,n)$ is not well understood\footnote{The Lusztig positive part of $\mcF(k,n)$ is well understood but it can be shown that the Pl\"ucker positive region $\mcF_+(k,n)$ strictly contains the Lusztig positive region when $n > 2k + 1$, see \cite{BBEG24}.} for general $k$. This motivates the following question:
     
     \begin{problem}
       Study the face structure of $\mcF_+(k,n)$ and find a parametrization of its cells for general $(k,n)$.
     \end{problem}

    The following result is known to the experts, see \cite[Proposition 5.1]{EPS24} and \cite[Equation (1.7)]{BL24}. We include a proof here for completeness.
    
	\begin{proposition}
		The homogeneous coordinate rings of the $2$-step flag variety $\mcF(k,2k+1)$ and the Grassmannian $\Gr(k+1, 2k+2)$ are isomorphic.
	\end{proposition}
	\begin{proof}

        Note that the cone over $\mcF(k,2k+1) \subset \PP^{\binom{2k+1}{k}-1} \times \PP^{\binom{2k+1}{k+1}-1}$ is an affine variety $\widehat{\mcF}(k, 2k+1) \subset \AA^{\binom{2k+1}{k}} \times \AA^{\binom{2k+1}{k+1}} = \AA^{\binom{2k+2}{k+1}}$ and the cone over the Grassmannian $\Gr(k+1, 2k+2)$ is the affine variety $\widehat{\Gr}(k+1, 2k+2) \subset \AA^{\binom{2k+2}{k+1}}$. To prove the above statement, we show that $\widehat{\mcF}(k,2k+1) = \widehat{\Gr}(k+1, 2k+2)$ in $\AA^{\binom{2k+2}{k+1}}$. To that end, it suffices exhibit a quasi affine variety $U$ that is open and dense in both $\widehat{\Gr}(k+1, 2k+2), \widehat{\mcF}(k,2k+1) \subset \AA^{\binom{2k+2}{k+1}}$. 
        
        The points of the variety $U \subset \AA^{\binom{2k+2}{k+1}}$ are obtained by taking taking the $(k+1) \times (k+1)$-minors of the $(k+1) \times (2k+2)$ matrices of the following form:
		\begin{equation}\label{eq:MatrixFlagVarGrass}
		    \begin{bmatrix}
			1 & \cdots & 0 & \ast & \cdots & \ast    & 0\\
			\vdots & \ddots & \vdots & \vdots & \cdots & \vdots & \vdots \\
			0 & \cdots & 1 & \ast &\cdots & \ast & 0\\
			0 & \cdots & 0 & {\color{blue} t} & \ast \cdots & \ast & {\color{red} z}
		\end{bmatrix},
		\end{equation}
		with $\ast \in \AA$ and $t,z \in \CC^{\times}$. It is clear how $U$ sits inside $\widehat{\Gr}(k+1, 2k+2)$. We now describe how to obtain points $(\widehat{V}, \widehat{W}) \in \widehat{\mcF}(k, 2k+1)$ from the matrices \eqref{eq:MatrixFlagVarGrass}. We obtain $\widehat{W}$ by taking $(k+1)\times(k+1)$ minors excluding the last column of \eqref{eq:MatrixFlagVarGrass}, while $\widehat{V}$ is obtained by taking the $k \times k$ minors of the upper left $k \times (2k+1)$ part of \eqref{eq:MatrixFlagVarGrass}, all multiplied by $z$. This shows that $\widehat{\mcF}(k,2k+1) = \widehat{\Gr}(k+1,2k+2)$ is the Zariski closure of $U$ in $\AA^{\binom{2k+2}{k+1}}$ which finishes the proof.
	\end{proof}
    
        \begin{remark}
            We warn the reader that although the homogeneous coordinate rings of $\mcF(k,2k+1)$ and $\Gr(k+1,2k+2)$ are isomorphic, the two varieties $\mcF(k,2k+1)$ and $\Gr(k+1,2k+2)$ are clearly \emph{not} isomorphic. Moreover, the coordinate ring of $\Gr(k+1,2k+2)$ has a $\ZZ$-grading while the coordinate ring of $\mcF(k,2k+1)$ has a $\ZZ^2$-grading.
        \end{remark}
        
	\begin{theorem}
		The positive orthogonal Grassmannians $\OGr_+(k,2k+1)$ and $\OGr_+(k+1,2k+2)$ can be identified through a linear isomorphism.
	\end{theorem}
	\begin{proof}
	We denote the Pl\"ucker coordinates on $\OGr(k,2k+1)$ by $p_I$ and those of $\OGr(k+1, 2k+2)$ by $q_J$. The orthogonal Grassmannian $\OGr(k,2k+1)$ is irreducible and the map:
    \begin{equation}\label{eq:Isom(k,2k+1)}
            \begin{matrix}
            \Phi_k \colon & \OGr(k+1, 2k+2) &   \to    & \OGr(k, 2k+1)  \\
                          & (q_J)_{J \in \binom{[2k+2]}{k+1}}           & \mapsto  &  (p_I = q_{I \cup \{2k+2\}})_{I \in \binom{[2k+1]}{k}}
            \end{matrix}
    \end{equation}
    is an isomorphism between the $\OGr(k,2k+1)$ (which is irreducible) and the standard component in $\OGr(k+1,2k+2)$. It is not so difficult to see that this isomorphism restricts to an isomorphism between the positive orthogonal Grassmannians $\OGr_+(k,2k+1)$ and $\OGr_+(k+1,2k+2)$.
	\end{proof}
	
	\begin{remark}
		The equations that cut out $\OGr(k,2k+1)$ in $\Gr(k,2k+1)$ are all quadrics. So it is remarkable that we can still describe the face structure of $\OGr_+(k,2k+1)$ from our understanding of the face structure of $\OGr_+(k+1, 2k+2)$ which is obtained by taking a linear slice of $\Gr_+(k+1,2k+2)$!
	\end{remark}
	
\begin{figure}[ht]
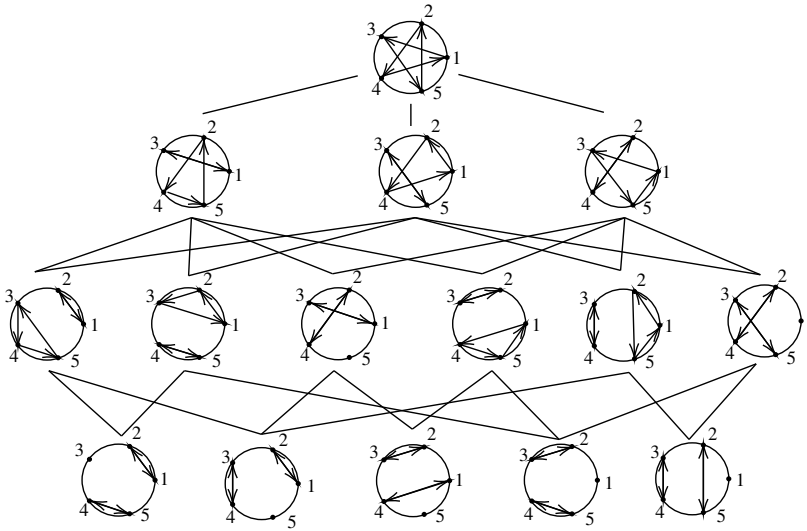

         \centering
        \scalebox{0.7}{
        \tikzset{every picture/.style={line width=0.75pt}} 

        }
         \caption{The face poset of $\OGr_+(2,5)$ matches that of $\OGr_+(3,6)$. See Figure 7 in \cite{galashin2020ising}. \label{fig:OGr25}}
     \end{figure}

	\begin{example}[$\OGr_+(2,5)$]
        The orthogonal Grassmannian $\OGr_+(2,5)$ is isomorphic to $\OGr_+(3,6)$. The Hasse diagram of the face poset of the latter is in \cite[Figure 7]{galashin2020ising}. Figure \ref{fig:OGr25} gives the same Hasse diagram in the realizable permutations in $\OGr_+(2,5)$. These cells can be parameterized using the isomorphism in \eqref{eq:Isom(k,2k+1)} and \cite[Theorem 5.17 (i)]{galashin2020ising}.
	\end{example}
 
     We finish this section by explaining how one goes from matchings $\tau$ on $[2k+2]$ to the admissible permutations in $[2k+1]$ i.e. permutations $\sigma$ of $[2k+1]$ with corresponding positroid cell $\Pi_{\sigma}$ such that $\Pi_{\sigma} \cap \OGr_+(k,2k+1)$ is nonempty.
     
\begin{figure}[ht]
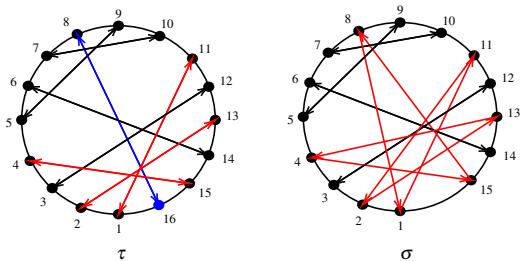

         \centering
         \scalebox{0.4}{    
    \tikzset{every picture/.style={line width=0.75pt}}

    }
         \caption{A matching $\tau$ of $[2k+2]$ and the corresponding permutation $\sigma$ of $[2k+1]$ for $k = 7$. In the left figure, starting vertex $16$ (in blue), the chords in red are longest sequence of chords $c_1,\dots, c_r$ that intersect each other pairwise. In the figure on the right, the blue vertex is deleted and the chords in red are turned into the unique cycle with support $\{1,2,4,8,11,13,15\}$ and 4 excedances.}
         \label{fig:matchingCorrespondence}
     \end{figure}

     Let $c$ denote the chord in $\tau$ attached to the vertex $2k+2$ and, starting from the vertex $2k+2$, consider the largest sequence $c = c_1, c_2, \dots, c_r$ of pairwise intersecting chords of $\tau$. Denote the $2r$ vertices of these chords by $i_1 < \dots < i_{2r-1} < 2k+2$. Then the cell $\Pi_{\tau} \cap \OGr_+(k+1, 2k+2)$ is isomorphic to the cell $\Pi_{\sigma} \cap \OGr_+(k, 2k+1)$ where $\sigma$  is the permutation of $[2k+1]$ obtained by replacing the chords $c_1, \dots, c_r$ with the unique cycle with support $\{i_1, \dots, i_{2r-1}\}$ and $r$ excedances. See Figure \ref{fig:matchingCorrespondence} for an example.
     
	\section{What goes wrong for \texorpdfstring{$\OGr_{+}(k, n)$}{} when \texorpdfstring{$n > 2k+1$}{} and \texorpdfstring{$k>1$}{}?}
	\label{sec5}
    In this section we show why positroid cells fail to induce a cell decomposition  of $\OGr_{+}(k, n)$ as soon as $n > 2k + 1$ and $k>1$. Let us start with the following:
    
	\begin{definition}\label{def:ortho-Positroids}
		For any positroid $\scrM$ of type $(k,n)$ and for any pair of subsets $I,J$ of $[n]$ of size $k-1$ we define the following two subsets of $[n]$:
		\[
		A_{IJ}^{\pm}(\scrM) = \big\{ \ell \in [n] \colon I\ell, J\ell \in \scrM \text{ and } (-1)^{\ell - 1} \epsilon_{I\ell} \epsilon_{J\ell} =  \pm 1 \big\}.
		\]
		We say that $\scrM$ is an orthopositroid if for any $I,J \in \binom{[n]}{k-1}$ we have:
		\[
		A_{IJ}^{+}(\scrM) = \emptyset \quad \iff \quad    A_{IJ}^{-}(\scrM) = \emptyset. 
		\]
	\end{definition}
 
	\begin{example}Let $n = 5$ and consider the two following positroids:
		\[
		\scrM_1 = \big\{\{1, 2\}, \{1, 4\}, \{2, 5\}, \{4, 5\} \big\} \ \text{and} \ \scrM_2 = \big\{ \{1, 2\}, \{1, 3\}, \{2, 4\}, \{3, 4\}\big\}.
		\]
		We then have:
		\[ 
		A^+_{24}(\scrM_1)  = \emptyset  \quad \text{and} \quad A^-_{24}(\scrM_1)  =  \{2\}.
		\]
		So $\scrM_1$ is \textbf{\emph{not}} an orthopositroid. One can check that $\scrM_2$ is an orthopositroid.
	\end{example}

    The motivation behind this definition is that if $X$ is a point in $\OGr_+(k,n)$ and $\mathscr{M}_X$ is its associated positroid then $\mathscr{M}_X$ is necessarily an orthopositroid in the sense of Definition \ref{def:ortho-Positroids}. This is because the Pl\"ucker coordinates of $X$ satisfy the equations \eqref{eq:orthRelationsAlt}.

    Since we will show that positroid cells dot not induce a CW cell decomposition of $\OGr_+(k,n)$, we refrain from elaborating more on the realizability of orthopositroids for general $k$, leaving that discussion for upcoming~work, and state the following

    \begin{conjecture}
        The orthopositroids in the sense of Definition \ref{def:ortho-Positroids} are realizable i.e. for any orthopositroid $\mathscr{M}$ of type $(k,n)$, there exists $X \in \OGr_+(k,n)$ such that $\mathscr{M} = \mathscr{M}_X$.        
    \end{conjecture}

    Let us start with $\OGr_+(2,6)$. A exhaustive computation shows that, out of all the positroids $\mathscr{M}$ (or decorated permutations $\sigma$) of type $(2,6)$, there are exactly $99$ orthopositroids (or admissible permutations). All of these orthopositroids are realizable in $\OGr_+(2,6)$ and are listed in Table \ref{tab:O26}. Let us now focus on the following orthopositroid cells $\sigma$ and $\tau$ in $\OGr_+(2,6)$:

    \begin{center}
    \scalebox{0.5}{
        \tikzset{every picture/.style={line width=0.75pt}} 
        

        }
    \end{center}
    
\begin{figure}[H]
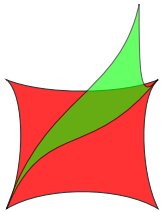

        \centering
        \begin{center}
        
				%
    \end{center}
        \caption{A cartoon of the cell $C_\sigma$ (in green) glued to the cell $C_\tau$ (in red).}
        \label{fig:gluing}
    \end{figure}

    We shall show that $C_\sigma = \Pi_{\sigma} \cap \OGr_+(2,6)$ and $C_\tau =  \Pi_{\tau} \cap \OGr_+(2,6)$, which are cells of dimension $2$, have the combinatorial type of a triangle and a square respectively. We start by giving generic matrices $M_\sigma, M_\tau$ that describe the points of the cells $C_\sigma, C_\tau$ respectively:
    \begin{align*}
        M_\sigma &= %
, \quad &\text{where } &\begin{cases}
                                  a,b,c>0 \\
                                  1 + b^2  = a^2 +  c^2\\
                              \end{cases}.
    \end{align*}
    Note that on the one hand, the closure of the cell $C_\tau$ in $\OGr_+(2,6)$ is isomorphic to $\OGr_+(1,4)$ so Theorem \ref{thm:OGr(1,n)} implies that $\overline{C_\tau}$ has the same combinatorial type as a square. On the other hand, the closure of the cell $C_\sigma$ has the combinatorial type of a triangle. Its edges are given by:
    \[
   \resizebox{.99\hsize}{!}{$
    e_1 = %
 \quad b \geq 0$}.
    \]
    The edge $e_1$ is one of the diagonals of the ``square" $C_\tau$. So the cell $C_\sigma$ glues with the cell $C_\tau$ as in Figure \ref{fig:gluing}.
    
   This shows that the positroid cells are not enough to induce a CW cell decomposition on $\OGr_+(2,6)$. In general this problem arises as soon as $n > 2k+1$.  This is because whenever $n > 2k+1$ we have $n-6 \geq  2(k-2)$, so we can extend a $2 \times 6$ matrix in $\OGr_+(2,6)$ by a $(k-2)\times(n-6)$ as follows
    \setcounter{MaxMatrixCols}{20}
    \[
        \resizebox{.7\hsize}{!}{$%
$.
    }
    \]
    
    We can then realize each positroid cell in $\OGr_+(2,6)$ as some positroid cell of $\OGr_+(k,n)$ and the same problem as above arises again. This highlights the need for new combinatorics to give a CW cell decomposition of $\OGr_+(k,n)$ when $n > 2k+1$ and $k > 1$.

    \begin{problem}
        Find a cell decomposition for $\OGr_+(k,n)$ when~$n > 2k+1$ and describe the combinatorics behind its face poset.
    \end{problem}

    \begin{table}[ht]
        \centering
        \scalebox{0.9}{
            \begin{tabular}{|c|c|c|}
                \hline
                dim. & permutations up to cyclic symmetry & \# \\
\hline
            \end{tabular}
            }
        \caption{The $99$ realizable decorated permutations in $\OGr_+(2,6)$. The dimensions of their cells are given in the leftmost column and the number of cells in each dimension in the rightmost column. \label{tab:O26}}
    \end{table}

    \footnotesize{{\bf Acknowlegements}: 
    We are grateful to Elizabeth Pratt for pointing us to the Weyl dimension formula \eqref{eq:wdf} used in the proof of Proposition \ref{prop:degree} and Theorem \ref{thm:Primality}. We also thank Lara Bossinger, Sebastian Seemann and the anonymous referee for helpful comments. YEM was partially supported by Deutsche Forschungsgemeinschaft SFB-TRR 195 “Symbolic Tools in Mathematics and their Application.” YM is supported by the National Science Foundation under Award DMS-2402069 and the Bob Moses Fund of the Institute for Advanced Study.}

    \bibliographystyle{acm}
    \bibliography{paper}
    
\end{document}